\documentclass[review,3p]{elsarticle}						

\usepackage{style_latex_share_writing}

\usepackage{style_latex_share_articles}

\usepackage{style_latex_amr}

\begin{document}

\title{On the best accuracy using the $h$-adaptive finite element refinement}
\author{Jie Liu}

\iftrue
\begin{abstract}

In \cite{liu2022practical}, a general algorithm is developed to efficiently obtain the best accuracy using the regular refinement.
The adaptive refinement allows for obtaining an accuracy with a smaller number of DoFs compared with the regular refinement.
In this paper, we investigate the best accuracy when using the adaptive refinement.
To this end, we study the evolution of the truncation error and the round-off error using the adaptive refinement.
For the former, a new threshold for the selection of the number of elements to be refined is proposed.
For the latter, the round-off error is quantified using the method proposed in \cite{liu2022practical}.
Moreover, for achieving a tolerance, we propose to use the line of the round-off error as a stopping criterion.


\end{abstract}
\fi

\iftrue
\begin{keyword}
Adaptive Finite Element Method, Order of Convergence, Optimal Thresholding, Efficiency.
\end{keyword}

\fi

  \maketitle

\section{Introduction}                  \label{section_introduction}

Many problems in applied sciences are modeled using (systems of) partial differential equation(s).
These equations are typically solved using numerical methods, such as finite difference (FD), finite volume (FV), and finite element (FE) methods.
In this paper, we will focus on the FE method, denoted by FEM.

\iftrue

The accuracy of the numerical approximation of the solution to the PDE is influenced by various factors~\cite{ferziger2012computational}.
Here, we will focus on two of the factors that depend opposite to the number of degrees of freedom, denoted by DoFs.
The first factor is related to the truncation error, which (theoretically) decreases with an increase in the number of DoFs. 
The second source of errors we consider is due to the round-off errors made due to the finite computer precision; 
this contribution typically increases with an increasing number of DoFs.

\fi

For a PDE problem on a domain, to achieve a tolerance, the computational mesh should be fine enough.
When the tolerance is not satisfied on the initial mesh, to improve the accuracy, one typically use $h$-refinement (reducing element size), $p$-refinement (increasing element order), or $hp$-refinement (combining the above methods).
In this paper, we focus on the $h$-refinement.








It contains two types: one is the regular refinement, denoted by REG, for which all elements are refined with each refinement;
the other is the adaptive refinement, denoted by AMR, for which only the elements of which the estimator is relatively large are refined, where a prescribed percentage is involved, denoted by pct~\cite{berger1984adaptive}.

The REG strategy quickly becomes impractical in 2D and, especially, 3D problems of practical interest as problem sizes quickly grow to millions of even billions of DoFs.
That is when the AMR strategy comes in handy to achieve a good resolution of local solution features while keeping the total number of DoFs at a manageable level.
\iftrue
However, the AMR strategy has the same difficulties with the finite computer precision as REG, that is, a locally too fine mesh can give rise to a significant round-off error that might, in the limit, dominate the overall error~\cite{alvarez2012round,liu2022practical}.
This paper will investigate this aspect of the AMR strategy that is hardly discussed in the literature.
\fi

The round-off error in the AMR strategy may be subject to the percentage used, and hence one should be careful with it.
Several approaches are pertinent to this problem.
In~\cite{dorfler1996convergent}, the percentage is chosen based on the reduction of the error estimator.                
The accumulation function of the estimator is analyzed in~\cite{dannenhoffer1987grid,pons2019adaptive}.
However, it is unclear which percentage is more friendly for reducing the round-off error.
In view of the above, the main aim of the paper is to design an AMR strategy that, given order $p$ of the elements, minimizes the total error by reducing the truncation error as much as possible without a strong increase in the round-off error. 






The paper is organized as follows.
In Section~\ref{section_model_problem_adaptive_refinement_ect}, the model problem and the FEM formulation are discussed.
In Section~\ref{section_new_requirements_on_pct}, for the AMR strategy, we introduce a novel strategy to find the optimal percentage and discuss the associated algorithm.
In Section~\ref{section_round_off_error_in_amr}, we investigate the round-off error using the AMR strategy.
Conclusions are drawn in Section~\ref{section_conclusion}.

\section{Model problem and FEM formulation}       \label{section_model_problem_adaptive_refinement_ect}

In Section~\ref{section_model_problem}, the model problems are illustrated.
In Section~\ref{section_basics_fem}, the FEM formulation is shown.



\subsection{Model problem}                  \label{section_model_problem}



\begin{figure}[!ht]
\centering
   \subfloat[Unit square\label{2d_geometry_shape_0}]{\includegraphics[width=0.2\linewidth]{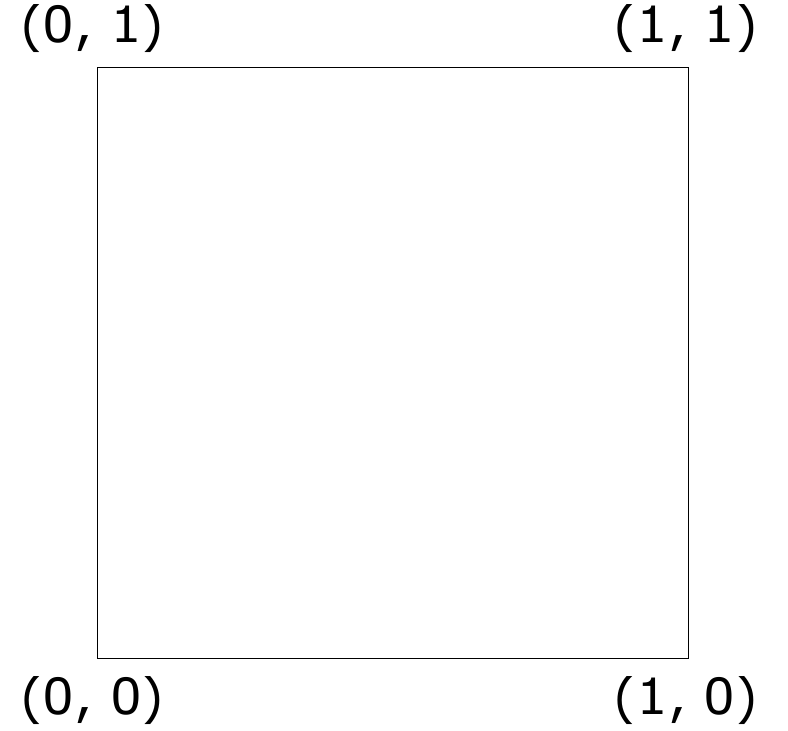}}  \quad
\caption{Shapes of the 2D geometry.}
\label{2d_geometry_shape}
\end{figure}

We consider one unit square, see Fig.~\ref{2d_geometry_shape}.
Dirichlet boundary conditions are imposed on the top and the bottom;
Neumann boundary conditions on the left and the right.


If not stated otherwise, we consider the Poisson problem with the following solution:

\begin{equation}
  u = e^{-\frac{(x-0.5)^2 + (y-0.5)^2}{c}},  \label{equation_2d_general}
\end{equation}
where $c$ is a constant.
When $c$ is equal to 1, the solution shape is very flat;
with the increase of $c$, the shape becomes sharper.
The shapes for $c = 1$ and $c$ = 1e-5 are shown in Fig.~\ref{shape_function_3d_solution_for_validation}.

\iftrue
\begin{figure}[!ht]
\centering
  \includegraphics[width=0.5\linewidth]{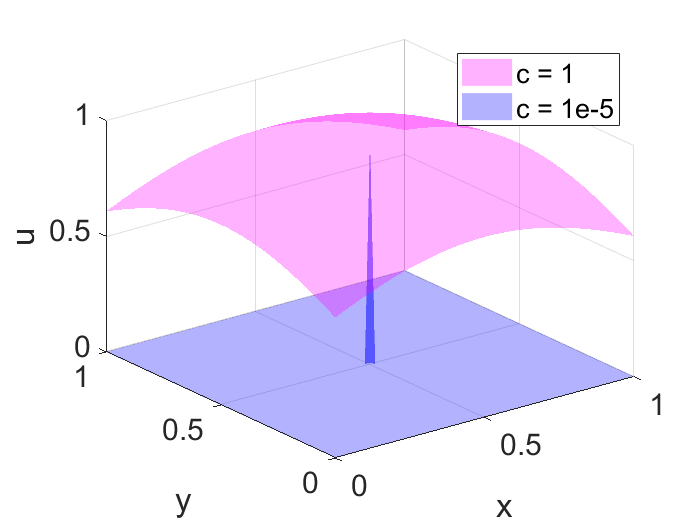}
  \caption{Shapes of the 2D solution on the unit square.}
  \label{shape_function_3d_solution_for_validation}
\end{figure}
\fi

\subsection{FEM formulation}              \label{section_basics_fem}

For the weak form of FEM, we refer to~\cite{liu2022practical}.
The function spaces are defined based on Lagrangian polynomials.
The element degree is denoted by $p$.


We implement FEM in deal.\rom{2}~\cite{dealII930}, where the IEEE-754 double precision is adopted.
It basically contains the following steps: generating the grid, setting up the problem, assembling the system of equations, solving, and computing the error.

For the details of generating the grid, see Section~\ref{section_basics_grid_generation}.
In the assembling process, for computing the occurring integrals, sufficiently accurate Gaussian quadrature formulas are used.
In the solving process, the UMFPACK solver~\cite{davis2004algorithm}, which implements the multi-frontal LU factorization approach, is used.
Using this solver prevents iteration errors associated with iterative solvers.
For the error estimation, we refer to Section~\ref{section_error_estimation_reg}.


The above steps involve the following parameters associated with a computational mesh:
the refinement level, grid size, number of DoFs, and error, which are denoted by $R$, $h$, $N$, and $E$, respectively.
$R = 0$ for the initial mesh, and it increases by one with each refinement.
Note that, $E$ can also be the truncation error $E _{\rm T}$ or the round-off error $E _{\rm R}$.

Moreover, the use of two types of refinement strategies ushers in two sets of parameters.
To differentiate the parameters of the adaptive refinement and the regular refinement, the notations in Table~\ref{table_differentiation_of_the_parameters_of_the_REG_and_the_AMR} are introduced.

\begin{table}[!ht]
\centering
\scriptsize
\caption{Differentiation of the parameters of the REG strategy and the AMR strategy.}
\begin{tabular}{c|c|c} \hline
 \diagbox{Parameter}{Refinement\\ strategy} & REG & AMR \\ \hline
 $R$ & $R ^{\rm reg}$ & $R ^{\rm amr}$ \\ \hline                            
 $h$ & $h ^{\rm reg}$ & $h ^{\rm amr}$ \\ \hline
 $N$ & $N ^{\rm reg}$ & $N ^{\rm amr}$ \\ \hline
 $E$ & $E ^{\rm reg}$ & $E ^{\rm amr}$ \\ \hline
\end{tabular}
\label{table_differentiation_of_the_parameters_of_the_REG_and_the_AMR}
\end{table}

\subsubsection{Grid generation}                 \label{section_basics_grid_generation}

For the adaptive refinement, the error estimator is chosen to be the Kelly estimator~\cite{kelly1983posteriori}, which is based on the gradient derived from the numerical solution.
For the mesh where we will conduct the adaptive refinement, we assume that the distribution of the error estimator is nonuniform~\cite{dorfler1996convergent}.
Note that the grid sizes of different cells are not the same since not all cells are refined with one refinement, and we use the minimum grid size to represent the overall grid size.


\subsubsection{Error estimation}                       \label{section_error_estimation_reg}

The numerical error depends on $h$ and $p$.
We can estimate the error either when $u_h$ is known or not, which are called a priori error estimation and a posteriori error estimation, respectively.


Using the former, the exact value of the error is not given, but an upper bound in provided.
It follows that~\cite{gockenbach2006understanding,chen2005finite}
\begin{equation}
  \| u - u_h \| _{L ^2} \leq C h ^{p+1} \| u \| _{H ^{p+1}},                    \label{formula_bound_of_the_truncation_error}
\end{equation}
where $C$ is a constant, and $u_h$ denotes the numerical solution.

\iftrue

Using the latter, the exact value of the error is computed, where the reference solution is the exact solution or that from other sources.
With respect to the latter, which is the general case, various methods, such as recovery methods, residual methods, and duality methods, can be used.
In this paper, we use the recovery method.
Specifically, the Richardson extrapolation where the reference solution is the finer solution with grid size $h/2$ is considered~\cite{chamoin2021pedagogical}.

\fi

\section{The optimal mesh for achieving a tolerance}                      \label{section_new_requirements_on_pct}



In Section~\ref{section_error_reduction_using_amr}, we illustrate the influence of pct on both the truncation error and the round-off error.
In Section~\ref{section_procedure_for_EOAMR}, we propose a new method to select pct.
For different problems, the numerical results of the new method can be found in Section~\ref{section_numerical_results_using_EOAMR}.

\subsection{Error reduction using the AMR strategy}              \label{section_error_reduction_using_amr}

\begin{figure}[!ht]
\centering
   \includegraphics[width=0.5\linewidth]{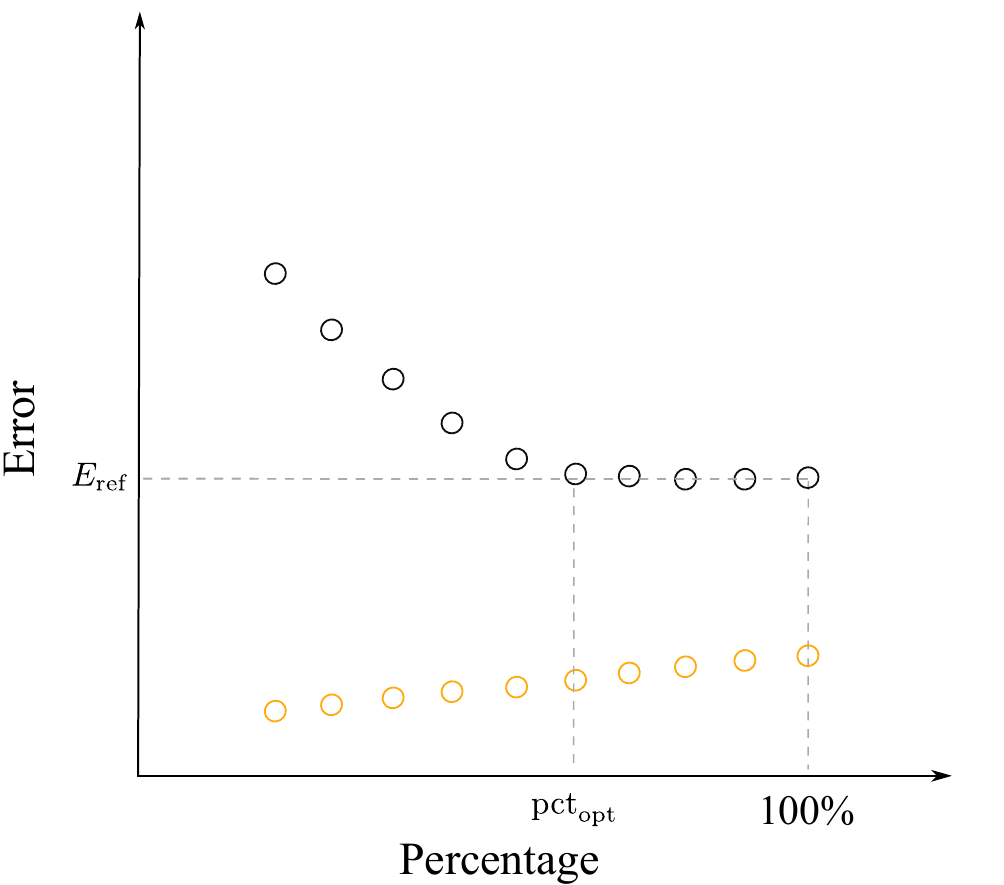}
   \caption{Influence of the pct on the error.}
   \label{sketch_influence_of_pct_on_error}
\end{figure}

When using the AMR strategy, different values of pct mean different degrees of error reduction.
Specifically, the resulting truncation error basically decreases with an increasing percentage;
if the pct is large enough, $E ^{\rm amr}$ becomes stable, see the black circles in Fig.~\ref{sketch_influence_of_pct_on_error}.
Obviously, the stabilized $E ^{\rm amr}$ is equal to $E ^{\rm reg}$.
Here, we denote the stabilized $E ^{\rm amr}$ by $E_{\rm ref}$, and the smallest pct that can achieve it by pct$_{\rm opt}$.

The pct$_{\rm opt}$ strongly depends on the distribution of the estimator.
If the distribution consists of very localized structures, it is very small, whereas for the distribution with a global structure, it is close to 100\%~\cite{pons2019adaptive}.

On the other hand, according to the results in Section~\ref{section_round_off_error_in_amr}, the round-off error increases with an increasing pct.
The sketch of their relation are denoted by orange circles in Fig.~\ref{sketch_influence_of_pct_on_error}.

\subsection{A new method for seeking pct}               \label{section_procedure_for_EOAMR}



For achieving a tolerance, to reduce the round-off error as much as possible, we propose to use pct$_{\rm opt}$.
As a first step, we obtain $E_{\rm ref}$ by conducting an extra regular refinement.
Second, starting with one small pct, we test the pct in ascending order until pct$_{\rm opt}$ is found, where the increment of pct for each test is chosen to be 10\% by experience.

For simplicity, the following abbreviations are introduced.
We denote the present mesh by $M_0$;
the pct used for each testing is denoted by pct$_{\rm tst}$, and the resulting error by $E_{\rm tst}$.
Specially, the initial pct$_{\rm tst}$ is denoted by pct$_{\rm init}$, and the resulting error by $E_{\rm init}$.
The selection of the pct$_{\rm init}$, and the procedure for testing different pct$_{\rm tst}$ can be found below.



\paragraph{Step 1: selecting \texorpdfstring{\rm pct$_{\rm init}$}{pctinit}}                  

It is slightly different for the first pct$_{\rm opt}$ and the remaining pct$_{\rm opt}$.
For the former, to capture the mesh locality as much as possible, pct$_{\rm opt}$ is chosen to be a relatively small number reading 5\%.
For the latter, we will choose the value based on the previous pct$_{\rm opt}$, denoted by pct$_{\rm prev}$.
That is, pct$_{\rm init}$ is chosen to be slightly smaller than pct$_{\rm prev}$ since the distribution of the error estimator tends to become flatter with each adaptive refinement, where the formula used reads:
\begin{equation}
  \text{pct}_{\rm init} = 0.7 \times \text{pct}_{\rm prev}.
\end{equation}

Note that if the resulting pct$_{\rm opt}$ is larger than 85\%, which indicates the distribution of the error estimator is smooth, we move to the REG strategy for the next refinement.

\paragraph{Step 2: seeking \texorpdfstring{\rm pct$_{\rm opt}$}{pctopt}}                  

When a pct$_{\rm init}$ is obtained, we can start searching for pct$_{\rm opt}$ by setting pct$_{\rm tst}$ = pct$_{\rm init}$.                

\begin{figure}[!ht]
\centering
   \includegraphics[width=0.5\linewidth]{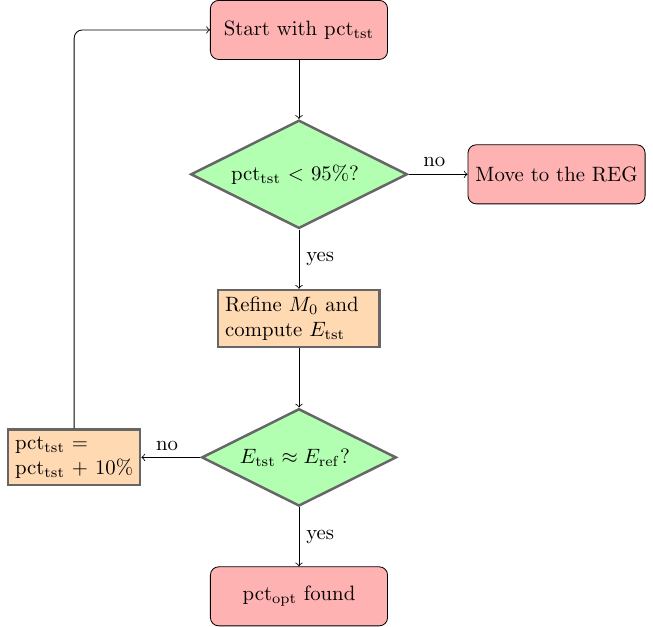}
   \caption{Illustration of the procedure for seeking pct$_{\rm opt}$.}
   \label{sketch_seeking_pct_opt}
\end{figure}

First, we refine $M_0$ using pct$_{\rm tst}$ and compute $E_{\rm tst}$.
Second, we compare $E_{\rm tst}$ with $E_{\rm ref}$:
if $E_{\rm tst} \approx E_{\rm ref}$, we are satisfied with this pct$_{\rm tst}$;
otherwise, we increase the pct$_{\rm tst}$.
With this new pct$_{\rm tst}$, if pct$_{\rm tst} < 95\%$, we repeat the above two stages;            
otherwise, we move to the regular refinement.
The associated algorithm can be found in Algorithm \ref{algo_selecting_pct}, also in Fig.~\ref{sketch_seeking_pct_opt}.

\vspace{0.2cm}
\begin{algorithm}[H]
\caption{Seeking pct$_{\rm opt}$}
\label{algo_selecting_pct}
pct$_{\rm tst}$ $\gets$ pct$_{\rm init}$\;                      
Refine $M_0$ and compute $E_{\rm tst}$\;
\While{$E_{\rm tst} \not\approx E_{\rm ref}$}
{
    pct$_{\rm tst}$ $\gets$ pct$_{\rm tst}$ + 10\%\;
    \eIf{\text{pct}$_{\rm tst} < 95\%$}
    {
        Refine $M_0$ and compute $E_{\rm tst}$\;
    }
    {
        Move to the regular refinement\;
    }
}
\end{algorithm}

In summary, we obtain pct$_{\rm opt}$ by adding extra work.                     
In what follows, we denote the $E_{\rm ref}$-oriented adaptive mesh refinement proposed above by EOAMR.

\subsection{Results}                \label{section_numerical_results_using_EOAMR}

\begin{figure}[!ht]
\centering
   \includegraphics[width=0.2\linewidth]{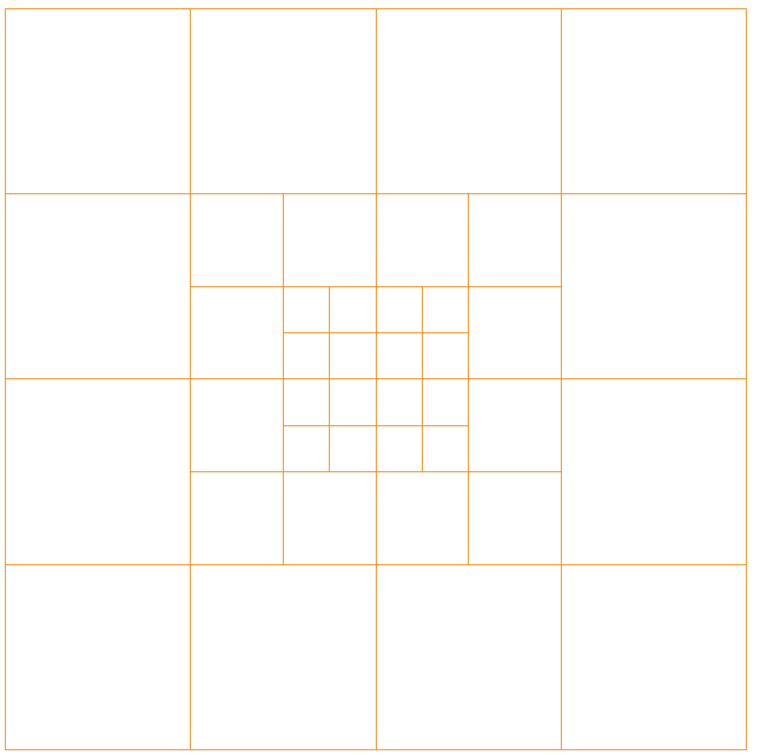}
   \caption{Initial mesh for the problem of Eq.~(\ref{equation_2d_general}) with $c=1e$-5.}
   \label{initial_mesh_fitted_to_the_load}
\end{figure}

We investigate the problem of Eq.~(\ref{equation_2d_general}) with $c$ = 1e-5.
The initial mesh that suits for the load can be found in Fig.~\ref{initial_mesh_fitted_to_the_load}.
As can be seen, since the load is very localized at the center, a smaller grid size is needed there.
The element degree $p$ used reads 1 and 3.
The results using the EOAMR can be found below.

\begin{figure}[!ht]
\centering
   \subfloat[$p=1$\label{comparison_of_EOAMR_to_REG_1p5_u_pct_deg_1}]{\includegraphics[width=0.33\linewidth]{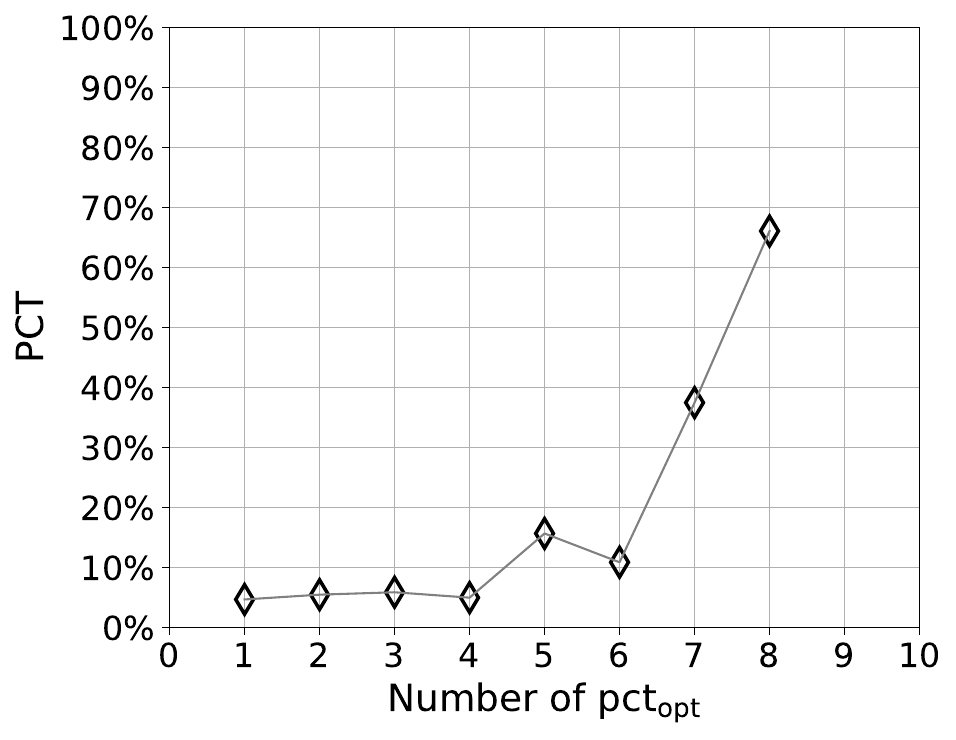}}
   \subfloat[$p=3$\label{comparison_of_EOAMR_to_REG_1p5_u_pct_deg_3}]{\includegraphics[width=0.33\linewidth]{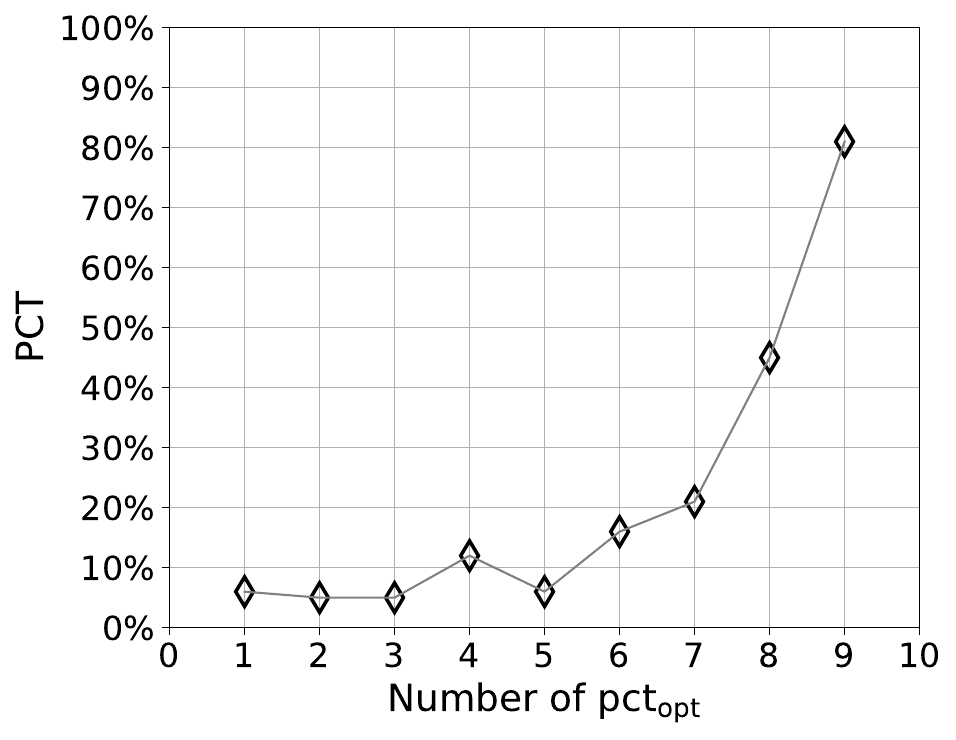}}
   \caption{Evolution of pct$_{\rm opt}$ using EOAMR for the 2D Poisson problem of Eq.~(\ref{equation_2d_general}) with $c=$1e-5.}
   \label{comparison_of_EOAMR_to_REG_1p5_u_pct}
\end{figure}

For $p = 1$, using the EOAMR, eight pct$_{\rm opt}$ are found.
The first pct$_{\rm opt}$ is much smaller than 1, reading about 5\%;
from the second one to the sixth one, pct$_{\rm opt}$ does not change much, which remains smaller than 20\%;
from the seventh one, pct$_{\rm opt}$ quickly increases, see Fig.~\ref{comparison_of_EOAMR_to_REG_1p5_u_pct_deg_1} for the evolution of pct$_{\rm opt}$.


%
%
For $p = 3$, nine pct$_{\rm opt}$ are found.
The evolution of pct$_{\rm opt}$ is similar with that for $p = 1$, see Fig.~\ref{comparison_of_EOAMR_to_REG_1p5_u_pct_deg_3}.


\sloppy

\section{Round-off error in \texorpdfstring{$h$}{h}-adaptive refinement}                   \label{section_round_off_error_in_amr}



We use the method of manufactured solutions proposed in~\cite{liu2022practical} to quantify the round-off error.
We consider the 2D Poisson problem with $u = 1$, and $p=2$ is investigated.
Two scenarios are analyzed. 
One is using different pct on the same mesh to see its influence on the error.
The other is using one constant pct to continually refine a mesh to study the error evolution. 
The details can be found below.

\paragraph{Using different pct on the same mesh}

\begin{figure}[!ht]
\centering
  \includegraphics[width=0.6\linewidth]{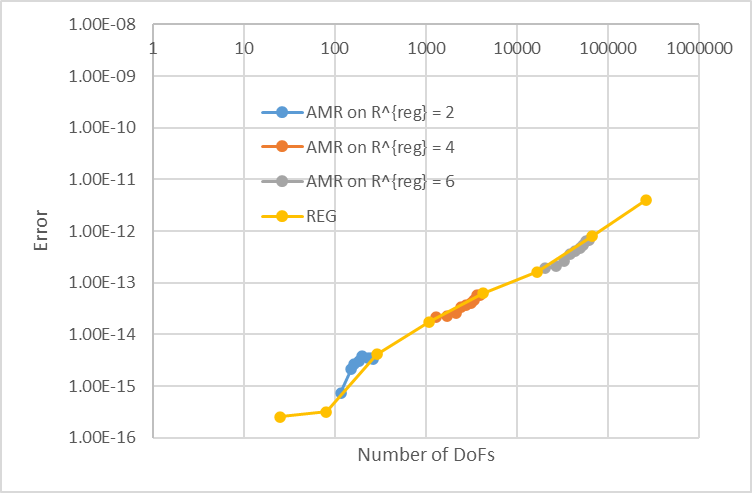}
  \caption{Errors using different pct on the same mesh.}
  \label{error_using_different_pct_on_one_refinement_level}
\end{figure}

We consider different meshes in terms of the regular refinement level $R^{\rm reg}$, of which the value equals 2, 4, and 6, respectively.
The pct ranges from 5\% to 85\%.
For each mesh, the error evolution with the number of DoFs using different pct can be found in Fig.~\ref{error_using_different_pct_on_one_refinement_level}.

As can be seen, for all the meshes, the error increases with the increasing pct.
When the pct is close to 1, the error nearly reaches that using the regular refinement.

\paragraph{Using one constant pct to refine one mesh continually}

\begin{figure}[!ht]
\centering
   \subfloat[1D \label{error_evolution_h_reg_vs_h_amr_p_equal_to_2_h_as_x_axis}]{\includegraphics[width=0.5\linewidth]{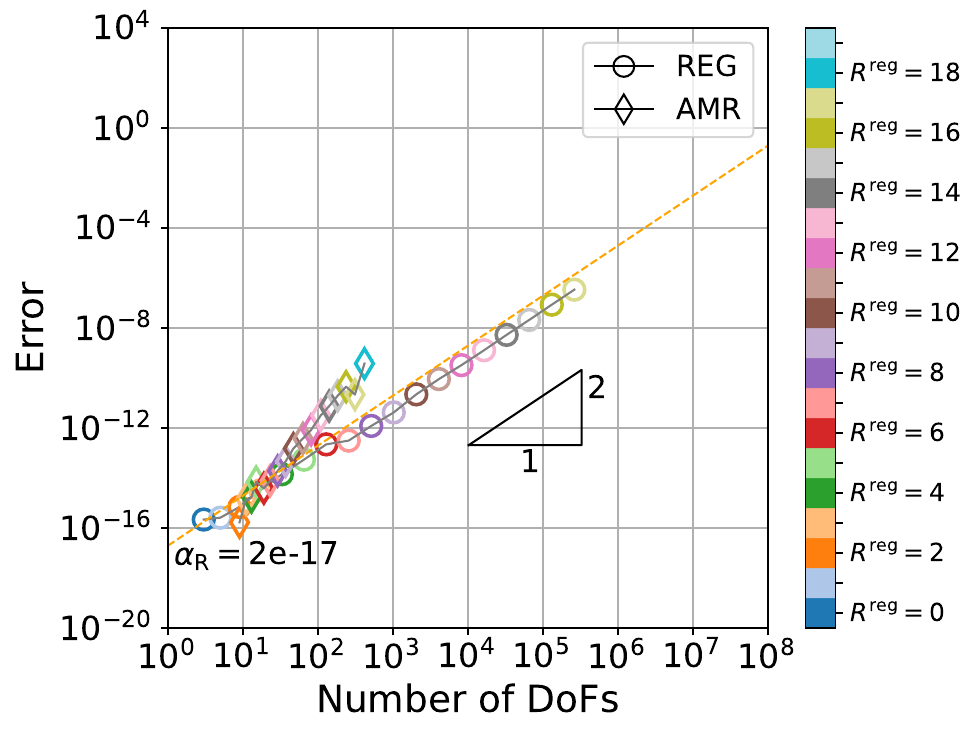}}
   \subfloat[2D \label{error_evolution_h_reg_vs_h_amr_p_equal_to_2_N_as_x_axis}]{\includegraphics[width=0.5\linewidth]{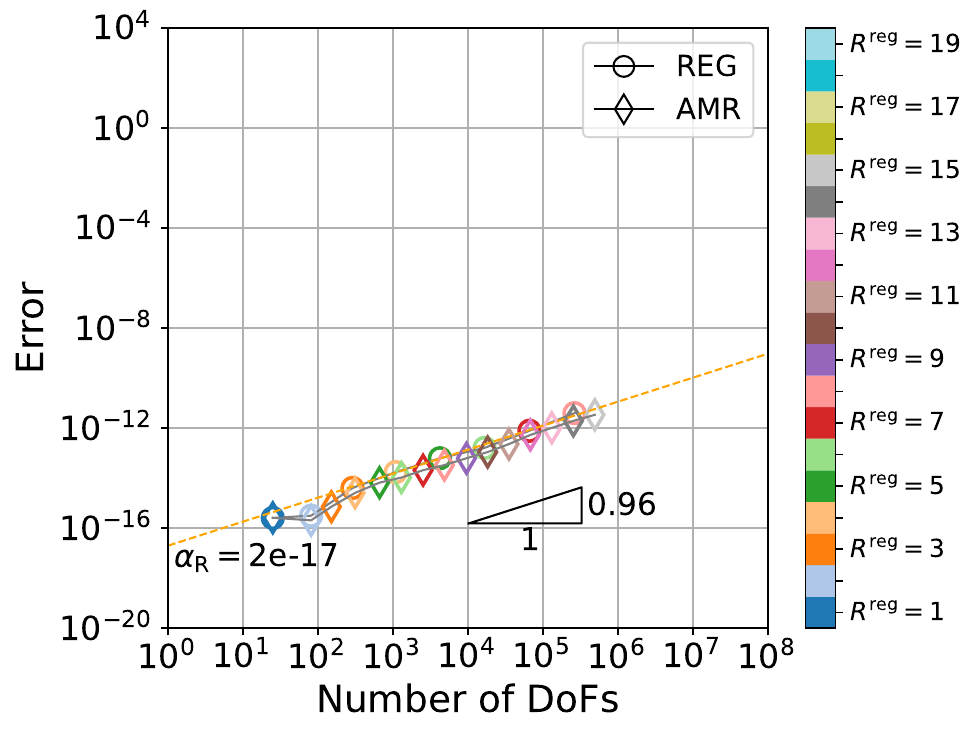}}
   \caption{Comparison of the error evolution using the regular refinement and the adaptive refinement.}
   \label{error_evolution_h_reg_vs_h_amr_p_equal_to_2}
\end{figure}

Here, pct = 30\% is chosen.
The initial mesh is the regular mesh of $R^{\rm reg} = 2$.
The resulting errors can be found in Fig.~\ref{error_evolution_h_reg_vs_h_amr_p_equal_to_2}.

Fig.~\ref{error_evolution_h_reg_vs_h_amr_p_equal_to_2_h_as_x_axis} indicates that by decreasing the grid size, $E _{\rm R}$ can either increase or decrease.
For the latter, it can be verified by comparing the error of the fifth regular mesh to that of the seventh adaptive mesh, both counting from the right.

From Fig.~\ref{error_evolution_h_reg_vs_h_amr_p_equal_to_2_N_as_x_axis}, $E _{\rm R}$ can be expressed by one function of $N$ when using the AMR strategy.
Obviously, the expression of $E _{\rm R}$ as a function of $N$ is basically the same for both the REG strategy and the AMR strategy.
However, this is not the case for 1D problems.                  



\section{Conclusion}                \label{section_conclusion}

We propose a new marking strategy that selects elements for refinement so that AMR follows the same asymptotic error behavior as REG without the excessive amount of DOFs of the latter. 
When the solution is sharp, using the new strategy to seek a tolerance, the CPU time can be saved a lot compared with using the regular refinement.
Our numerical experiments furthermore demonstrate that the round-off error does not depend on the minimum mesh width but on the total number of DoFs (i.e. it is a result of how many times we commit a round-off error). Hence, AMR not only reduces the amount of DOFs in regions where they are not needed to improve the accuracy but also mitigates the devastating effect of the round-off error and thereby enables more accurate solutions.
This advantage of AMR becomes particularly pronounced for high-order approximations.





\appendix


\section{A practical algorithm for finding \texorpdfstring{$h_{\rm opt}$}{hopt} when using the REG strategy}         \label{section_practical_alg_for_reg}

Normally, to reach one tolerance, we check if the tolerance is satisfied for each refinement level and stop further refinements when the tolerance is satisfied.
Since $tol$ is reachable when $tol \geq E_{\rm min}$, one may waste a lot of CPU time if $tol$ cannot be reached.

To avoid unnecessary computations, we predict $N_{\rm opt}$ and the associated $E_{\rm min}$ in \cite{liu2022practical}, of which the formula is derived from the perspective of the round-off error.
Here, we will illustrate a slightly different version: to predict $h_{\rm opt}$ and the associated $E_{\rm min}$ from the perspective of satisfying the condition of the truncation error, which is more concise and easy to understand.
For this strategy, the preliminaries and procedure can be found in \ref{section_preliminaries_for_prediction_reg} and \ref{section_procedure_for_prediction_reg}, respectively.
We show the application of the algorithm in Section~\ref{section_application_of_the_algorithm}.

\subsection{Preliminaries}                    \label{section_preliminaries_for_prediction_reg}




We use the following properties of $E_{\rm T}$ and $E_{\rm R}$.

\paragraph{Truncation error}

The a priori error estimation provides the theoretical convergence order of the error, which is for the validation purpose.
Specifically, when $h$ is relatively small, the error converges to a formula~\cite{convergenceFEM20xxscadoffice,gockenbach2006understanding}:
\begin{equation}
  E_{\rm T} = C_{\rm T} \cdot h^{q},                      \label{formular_of_E_T_in_terms_of_h}
\end{equation}
where $C_{\rm T}$ is a constant dependent on the element degree, and $q$ the analytical order of convergence, which equals $p+1$ when $p$ is smaller than ten approximately~\cite{mitchell2015high}.

\paragraph{Round-off error}             \label{section_error_estimation_reg_E_R}

\iftrue

The round-off error is caused by the adoption of the finite precision arithmetic, of which the value is not reflected in Eq.~(\ref{formula_bound_of_the_truncation_error}).
According to~\cite{liu386balancing}, it has a power-law relation with the number of DoFs.
In the log-log plot where the $x$ axis is the number of DoFs, and the $y$ axis is the error, the round-off error increases along a straight line with the increasing number of DoFs.
Denoting the slope and the offset of the line by $\alpha_{\rm R}$ and $\beta _{\rm R}$, respectively, the round-off error reads
\begin{equation}
 E_{\rm R} = \alpha_{\rm R} \cdot {N}^{\beta_{\rm R}}.              \label{formula_round_off_error_in_terms_of_N}
\end{equation}
In \cite{liu386balancing}, we show that $\beta_{\rm R}$ is independent of different $f$, and $\alpha_{\rm R}$ is linearly proportional to $\| u \|_2$.
For general problems, $\alpha_{\rm R}$ and $\beta_{\rm R}$ can be determined using the method of manufactured solutions~\cite{liu2022practical}.

\fi

\iftrue
Moreover, when $h$ is relatively small, $N$ is a function of $h$ and $p$~\cite{liu2022practical}:
\begin{equation}
  N \approx \Bigl(\frac{p}{h}\Bigr)^2,                     \label{formula_of_N_in_terms_of_h}
\end{equation}
which indicates that the number of DoFs basically quadruples with the increasing refinement level~\cite{liu2022practical}.
\fi
Therefore, substituting Eq.~(\ref{formula_of_N_in_terms_of_h}) into Eq.~(\ref{formula_round_off_error_in_terms_of_N}), we have the expression of $E_{\rm R}$ in terms of $h$:
\begin{equation}
 E_{\rm R} = C_{\rm R} \cdot {h}^{D_{\rm R}},              \label{formula_round_off_error_in_terms_of_h}
\end{equation}
where $C _{\rm R} = \alpha_{\rm R} \cdot p ^{2 \beta_{\rm R}}$, and $D_{\rm R} = -2 \beta_{\rm R}$.

Now, from Eq.~(\ref{formular_of_E_T_in_terms_of_h}) and Eq.~(\ref{formula_round_off_error_in_terms_of_h}), when $q$ is satisfied, the total error in terms of $h$ reads
\begin{equation}
  E = C_{\rm T} \cdot h^{q} + C_{\rm R} \cdot {h}^{D_{\rm R}}.                      \label{formular_of_E_in_terms_of_h}
\end{equation}
As can be seen, when $h$ is small, it basically follows that $E_{\rm T} > E_{\rm R}$;
with the decreasing $h$, $E_{\rm T}$ decreases, and $E_{\rm R}$ increases.

Due to the interplay between $E_{\rm T}$ and $E_{\rm R}$, there is a minimum achievable error, which is denoted by $E_{\rm min}$, for which the corresponding $h$ is denoted by $h_{\rm opt}$.

\subsection{Procedure for the prediction}                    \label{section_procedure_for_prediction_reg}


From Eq.~(\ref{formular_of_E_in_terms_of_h}), we can predict the error when order $q$ is satisfied.
First, we determine when $q$ is reached.
Second, we derive the value of $C_{\rm T}$.
Finally, we conduct the prediction.

\paragraph{Checking the condition of \texorpdfstring{$q$}{q}}                \label{section_criterion_for_the_prediction}

The practical value of $q$, denoted by $q_h$, is computed as follows.
\begin{equation}
  q_h = \frac{\log \left( \frac{E_{h/2}}{E_{h}} \right)}{\log 2},  \label{formula_order_of_convergence}
\end{equation}
where $E_h$ and $E_{h/2}$ denote the errors with grid size $h$ and $h/2$, respectively.

When order $q$ is reached, error $E$, number of DoFs $N$, refinement level $R$ and grid size $h$ are denoted by $E_{\rm c}$, $N_{\rm c}$, $R_{\rm c}$ and $h_{\rm c}$, respectively.
For different problems, $R_{\rm c}$ are different~\cite{Runborg2012VerifyingNC}.    
When the solution is smooth, $R_{\rm c}$ is relatively small.
Note that in this paper, we assume $q$ is reachable if not stated otherwise, and the CPU time required is not very large.


\paragraph{Deriving the coefficient}

The coefficient $C_{\rm T}$ is computed as follows.
\begin{equation}
  C_{\rm T} = \frac{E_{\rm c}}{{h_{\rm c}} ^q}.
\end{equation}

\paragraph{Predicting the results}                  \label{section_prediction_results}

For the tolerance $tol$, the corresponding $N$, $h$, and $R$ are denoted by $N_{tol}$, $h_{tol}$, and $R_{tol}$, respectively.
Specifically, when $tol = E_{\rm min}$, they are denoted by $N_{\rm opt}$, $h_{\rm opt}$, and $R_{\rm opt}$, respectively.

\begin{subequations}
For Eq.~(\ref{formular_of_E_in_terms_of_h}), $h_{\rm opt}$ and $E_{\rm min}$ are predicted as follows:
\begin{align}
 h_{\rm opt} ^{\rm PRED}&= \left( - \frac{C_{\rm T} \cdot q}{C _{\rm R} \cdot D_{\rm R}} \right)^{\frac{1}{D_{\rm R} - q}}, \label{formula_N_opt_in_terms_of_h}  \\
 E_{\rm min} ^{\rm PRED}&= C_{\rm T} \cdot {h_{\rm opt}} ^q + C_{\rm R} \cdot {h_{\rm opt}}^{D _{\rm R}},        \label{formula_E_min_in_terms_of_h} 
\end{align}
where $E_{\rm min} ^{\rm PRED}$ and $N_{tol} ^{\rm PRED}$ denote the predicted $E_{\rm min}$ and $N_{\rm opt}$, respectively.
\label{formula_h_opt_and_E_min_in_terms_of_h}%
Moreover, using the BF approach, $E_{\rm min}$ and $N_{\rm opt}$ are denoted by $E_{\rm min} ^{\rm BF}$ and $N_{tol} ^{\rm BF}$, respectively.
\end{subequations}

\subsection{Application}                   \label{section_application_of_the_algorithm}
                                                                    
Using the PRED method, we can relatively quickly determine if an accuracy can be reached.
For a given $tol$, one $p$ is possible when $tol \geq E_{\rm min} ^{\rm PRED}$.
With respect to $h_{tol}$, it is predicted by
\begin{equation}
  h_{tol} ^{\rm PRED} = {\Big( \frac{tol}{C_{\rm T}} \Big)}^{\frac{1}{q}},         \label{formular_h_tol_PRED_in_terms_of_h}
\end{equation}
where $h_{tol} ^{\rm PRED}$ denotes the predicted $h_{tol}$.
With $h_{tol} ^{\rm PRED}$, we can directly compute the solution of interest, and hence several refinement steps before $h_{tol} ^{\rm PRED}$ can be reduced compared to REG.
When $tol$ cannot be reached, we stop further refinements.

The reduction of the refinement steps indicates the reduction of the CPU time.
Specifically, when $tol$ can be reached, the CPU time reduction is not much compared with the BF method because we still have to compute the solution for $h_{tol} ^{\rm PRED}$, and the CPU time reduced is not as much as that for computing the solution for $h_{tol} ^{\rm PRED}$;
otherwise, the reduction of the CPU time is much more.

\bibliographystyle{unsrt}
\bibliography{bibfile}

\end{document}